\documentclass{amsart}
\usepackage{latexsym,amssymb,amsmath,amsthm,amscd,graphicx}

\usepackage{esint} 
\usepackage{color}

\numberwithin{equation}{section}
\newtheorem{theorem}{Theorem}[section]

\theoremstyle{definition}

\theoremstyle{remark}

\usepackage{url,hyperref}

\usepackage{color}
\definecolor{dblue}{rgb}{0,0,0.45}
\definecolor{red}{rgb}{0.7,0,0}

\begin{document}

\title[Three geometric constants for Morrey spaces]
{Three geometric constants for Morrey spaces}

\author[H.~Gunawan]{Hendra~Gunawan}
\address{Analysis and Geometry Group, Faculty of Mathematics
and Natural Sciences, Bandung Institute of Technology,
Bandung 40132, Indonesia.}
\email{hgunawan@math.itb.ac.id}

\author[E.~Kikianty]{Eder~Kikianty}
\address{Address: Department of Mathematics and Applied Mathematics,
University of Pretoria, Private Bag X20, Hatfield 0028, South Africa}
\email{eder.kikianty@gmail.com}

\author[Y.~Sawano]{Yoshihiro~Sawano}
\address{Department of Mathematics and Information Science,
Tokyo Metropolitan University, Minamioosawa, Hachouji-city 1-1,
Tokyo 192-0364, Japan
+
Department of Mathematics Analysis and the Theory of functions,
Peoples' Friendship University of Russia, Moscow, Russia
}
\email{ysawano@tmu.ac.jp}

\author[C.~Schwanke]{Christopher~Schwanke}
\address{Address: Department of Mathematics,
Lyon College, Batesville, AR 72501, USA
+
Unit for BMI, North-West University, Private Bag X6001,
Potchefstroom 2520, South Africa}
\email{cmschwanke26@gmail.com}

\subjclass[2010]{46B20}

\keywords{Morrey spaces, Dunkl-Williams constant, James constant, von Neumann-Jordan constant}

\begin{abstract}
In this paper we calculate three geometric constants, namely the von Neumann-Jordan constant,
the James constant, and the Dunkl-Williams constant, for Morrey spaces and discrete Morrey
spaces. These constants measure uniformly nonsquareness of the associated spaces. We obtain
that the three constants are the same as those for $L^1$ and $L^\infty$ spaces.
\end{abstract}

\maketitle

\section{Introduction}

The {\it von Neumann-Jordan constant} ${\rm C}_{\rm N J}(X)$ (see \cite{JoNe35}),
the {\it James constant} ${\rm C}_J(X)$ (see \cite{James64})
and
the {\it Dunkl-Williams constant} ${\rm C}_{\rm D W}(X)$ (see \cite{DuWi64})
for a Banach space $X$ are given by
$$
{\rm C}_{\rm N J}(X):=
\sup\left\{
\frac{\|x+y\|^2_X+\|x-y\|^2_X}{2(\|x\|^2_X{}+\|y\|^2_X)}
\,:\,x,y \in X \setminus \{0\}
\right\},
$$
$$
{\rm C}_{\rm J}(X):=
\sup\left\{\min \{\|x+y\|_X,\|x-y\|_X\}\,:\,
x,y \in X, \, \|x\|_X=\|y\|_X=1\right\},
$$
and
\[
{\rm C}_{\rm D W}(X):=
\sup\left\{
\frac{\|x\|_X+\|y\|_X}{\|x-y\|_X}
\left\|\frac{x}{\|x\|_X}-\frac{y}{\|y\|_X}\right\|_X\,:\,
x,y \in X,\, x \ne 0, y \ne 0, x \ne y
\right\},
\]
respectively. It is well known that $1\le C_{\rm N J}(X)\le 2$ for every Banach
space $X$, and that $C_{\rm N J}(X)=1$ if and only if $X$ is a Hilbert space.
Meanwhile, $\sqrt{2}\le C_{\rm J}(X)\le 2$ holds for every Banach space $X$,
and $C_{\rm J}(X)=\sqrt{2}$ if (but not only if) $X$ is a Hilbert space
(see \cite{Cassini86, GaLa90}).
As for the Dunkl-Williams constant, we have
$2\le {\rm C}_{\rm D W}(X) \le 4$  and ${\rm C}_{\rm D W}(X)=2$
if and only if $X$ is a Hilbert space \cite{DuWi64}.
For Lebesgue spaces $L^p=L^p({\mathbb R}^d)$ where $1\le p\le\infty$,
we have $C_{\rm N J}(L^p)=\max \{2^{2/p-1},2^{1-2/p}\}$ and
$C_{\rm J}(L^p)=\max \{2^{1/p},2^{1-1/p}\}$ \cite{KM01}.
Meanwhile, we know that $C_{\rm D W}(L^1)=C_{\rm D W}(L^\infty)=4$
\cite{JLM08}.

In this paper, we shall calculate the three constants for Morrey spaces
and discrete Morrey spaces. Let $1 \le p \le q<\infty$. The {\it Morrey space}
${\mathcal M}^p_q={\mathcal M}^p_q({\mathbb R}^d)$
is the set of all the measurable functions $f$ on ${\mathbb R}^d$ for which
$$
\| f \|_{{\mathcal M}^p_q}
:=
\sup_{B=B(a,r)}
|B|^{\frac1q-\frac1p}\left(\int_B |f(y)|^p\,dy\right)^{\frac1p}
<\infty,
$$
where $B(a,r)$ denotes the ball centered at $a\in{\mathbb R}^d$
having radius $r>0$ and Lebesgue measure $|B|$ (see, e.g., \cite{AX12}).
Since ${\mathcal M}^p_q$ is a Banach space, it follows from
\cite{Cassini86,DuWi64,GaLa90} that
\[
C_{\rm N J}({\mathcal M}^p_q),\ C_{\rm J}({\mathcal M}^p_q)\le 2\ \
{\rm and}\ \ C_{\rm D W}({\mathcal M}^p_q)\le 4.
\]

Our result for Morrey spaces is the following:

\begin{theorem}\label{thm1}
If $1 \le p<q<\infty$,
then $C_{\rm N J}({\mathcal M}^p_q)=C_{\rm J}({\mathcal M}^p_q)=2$
and $C_{\rm D W}({\mathcal M}^p_q)=4$.
\end{theorem}

Note that ${\mathcal M}^p_p=L^p$ holds and that their norms are identical.
The above theorem tells us that the case where $q>p$ is quite different
from the case where $q=p$. When $q>p$, the three constants
$C_{\rm J}({\mathcal M}^p_q)$, $C_{\rm N J}({\mathcal M}^p_q)$, and
$C_{\rm D W}({\mathcal M}^p_q)$ take the same value as those for $L^1$ and
$L^\infty$ spaces.

Moving on to discrete Morrey spaces, let $\omega:=\mathbb{N}\cup \{0\}$.
For $m:=(m_1,\dots, m_d)\in\mathbb{Z}^d$ and $N\in \omega$, let
$$
S_{m,N}:= \{k\in \mathbb{Z}^d\colon \|k-m\|_\infty \leq N\},
$$
where $\|(m_1, \dots, m_d)\|_\infty:=\max \{|m_i|\colon 1\leq i\leq d\}$ for
$(m_1, \dots, m_d)\in \mathbb{Z}^d$. The cardinality of $S_{m,N}$, denoted by
$|S_{m,N}|$, is $(2N+1)^d$, for every $m\in \mathbb{Z}^d $ and $N\in\omega$.
Given $1\leq p\leq q<\infty$, we define the {\it discrete Morrey space}
$\ell^p_q=\ell^p_q(\mathbb{Z}^d)$ to be the space of all functions
(sequences) $x\colon \mathbb{Z}^d
\rightarrow \mathbb{R}$ for which
$$
\|x\|_{\ell^p_q}:=\sup_{m\in\mathbb{Z}^d, N\in \omega} |S_{m,N}|^{\frac1q-\frac1p}
\left(\sum_{k\in S_{m,N}} |x(k)|^p\right)^\frac1p<\infty.
$$
We note that $\ell^p_q$, equipped with the above norm, is a Banach space (see
\cite{GKS18}). Our result for discrete Morrey spaces is the following:

\begin{theorem}\label{thm2}
If $1 \le p<q<\infty$, then $C_{\rm N J}(\ell^p_q)=C_{\rm J}(\ell^p_q)=2$
and $C_{\rm D W}(\ell^p_q)=4$.
\end{theorem}

This theorem also tells us that the case where $q>p$ is quite different from
the case where $q=p$ (where $\ell^p_p=\ell^p$).

\section{Proof of Theorems}

We prove both theorems by finding two elements in the space such that
the associated expressions are equal to two, two, and four, respectively.

\subsection{Proof of Theorem \ref{thm1}}

\begin{proof}
Let $1\leq p<q<\infty$, and let $f(x):=|x|^{-d/q}$, $x \in {\mathbb R}^d$,
where $|x|$ denotes the Euclidean norm of $x$. Then $f \in
{\mathcal M}^p_q({\mathbb R}^d)$ (see {\cite[\S 2]{SaSuTa11}}).
Define $g(x) := \chi_{(0,1)}(|x|)f(x)$,
$h(x) :=f(x)-g(x)$, and $k(x):=-f(x)+2g(x)$, for $x \in {\mathbb R}^d$.
By a change of variables, we see that
$$\|t^{d/q}g(t\cdot)\|_{{\mathcal M}^p_q}=\|g\|_{{\mathcal M}^p_q}$$
and
$$\|t^{d/q}h(t\cdot)\|_{{\mathcal M}^p_q}=\|h\|_{{\mathcal M}^p_q}$$
for all $t>0$. Since
\[
t^{d/q}g(t x)= \chi_{(0,1)}(t|x|)f(x)
\]
and
\[
t^{d/q}h(t x)= \chi_{(0,1)}(t|x|)f(x)-\chi_{[1,\infty)}(t|x|)f(x)
\]
for $t>0$ and $x \in {\mathbb R}^d$, by the monotone convergence property of
Morrey spaces we have
\[
\|f\|_{{\mathcal M}^p_q}=\|g\|_{{\mathcal M}^p_q}=\|h\|_{{\mathcal M}^p_q}=
\|k\|_{{\mathcal M}^p_q} \in (0,\infty).
\]
This implies that
$$\|f+k\|^2_{{\mathcal M}^p_q}{}+\|f-k\|^2_{{\mathcal M}^p_q}=
4(\|f\|^2_{{\mathcal M}^p_q}{}+\|k\|^2_{{\mathcal M}^p_q}{})$$
and
$$
\min \{\|f+k\|_{{\mathcal M}^p_q},\|f-k\|_{{\mathcal M}^p_q}\}
=
\min \{\|2g\|_{{\mathcal M}^p_q},\|2h\|_{{\mathcal M}^p_q}\}=
2\|f\|_{{\mathcal M}^p_q}=2\|k\|_{{\mathcal M}^p_q}.
$$
By definition and the fact that both $C_{\rm NJ}({\mathcal M}^p_q),\
C_{\rm J}({\mathcal M}^p_q)\le 2$,
we conclude that
\[
C_{\rm N J}({\mathcal M}^p_q)=C_{\rm J}({\mathcal M}^p_q)=2,
\]
as desired.

Finally, we calculate the Dunkl--Williams constant using the same ideas as in
\cite{JLM08}.
We consider $f$ and $(1+r)g+(1-r)h$ for $r \in (0,1)$.
We calculate
\begin{align*}
&\frac{\|f\|_{{\mathcal M}^p_q}+\|(1+r)g+(1-r)h\|_{{\mathcal M}^p_q}}{\|f-(1+r)g-(1-r)h\|_{{\mathcal M}^p_q}}
\left\|\frac{f}{\|f\|_{{\mathcal M}^p_q}}-\frac{(1+r)g+(1-r)h}{\|(1+r)g+(1-r)h\|_{{\mathcal M}^p_q}}\right\|_{{\mathcal M}^p_q}\\
&=\frac{\|f\|_{{\mathcal M}^p_q}+(1+r)\|f\|_{{\mathcal M}^p_q}}{r\|f\|_{{\mathcal M}^p_q}}
\left\|\frac{f}{\|f\|_{{\mathcal M}^p_q}}-\frac{(1+r)g+(1-r)h}{(1+r)\|f\|_{{\mathcal M}^p_q}}\right\|_{{\mathcal M}^p_q}\\
&=\frac{\|f\|_{{\mathcal M}^p_q}+(1+r)\|f\|_{{\mathcal M}^p_q}}{r\|f\|_{{\mathcal M}^p_q}}
\left\|\frac{2r h}{(1+r)\|f\|_{{\mathcal M}^p_q}}\right\|_{{\mathcal M}^p_q}\\
&=
\frac{4+2r}{1+r}.
\end{align*}
If we let $r \downarrow 0$, we obtain $C_{\rm D W}({\mathcal M}^p_q)=4$,
as required.
\end{proof}

Before we conclude this subsection, a remark may be in order.
Let $1 \le p \le q<\infty$. The {\it local Morrey space}
${\rm L}{\mathcal M}^p_q={\rm L}{\mathcal M}^p_q({\mathbb R}^d)$
is the set of all the measurable functions $f$ on ${\mathbb R}^d$
for which
$$
\| f \|_{{\rm L}{\mathcal M}^p_q}
:=
\sup_{B=B(0,r)}
|B|^{\frac1q-\frac1p}\left(\int_B |f(y)|^p\,dy\right)^{\frac1p}
<\infty.
$$
Arguing similary as before, we see that
$C_{\rm N J}({\rm L}{\mathcal M}^p_q)=C_{\rm J}({\rm L}{\mathcal M}^p_q)=2$
and $C_{\rm D W}({\rm L}{\mathcal M}^p_q)=4$
whenever $1 \le p<q<\infty$.

\subsection{Proof of Theorem \ref{thm2}}

\begin{proof}
Let $1\leq p<q<\infty$, and let us first consider the case where $d=1$.
Let $n\in\mathbb{Z}$ be an even
number with
$
n> 2^{\frac{q}{q-p}}-1,
$
or equivalently
$$
(n+1)^{\frac1q-\frac1p}<2^{-\frac1p}.
$$
Consider the sequence $(x_k)_{k\in \mathbb{Z}}$ defined by
$$
x_{0}=x_{n}=1, \text{and} \ x_k=0 \text{ for all } k\not\in \{0,n\}
$$
and the sequence $(y_k)_{k\in \mathbb{Z}}$ defined by
$$
y_{0}=1,\ y_{n}=-1,  \text{and} \ y_k=0 \text{ for all }k\not\in \{0,n\}.
$$
Then, we have
\begin{align*}
\nonumber \|x\|_{\ell^p_q}
&=\sup_{m\in\mathbb{Z}, N\in \omega} |S_{m,N}|^{\frac1q-\frac1p}
\left(\sum_{k\in S_{m,N}} |x_k|^p\right)^\frac1p\\
&=\max\left\{ 1, |S_{\frac{n}{2},\frac{n}{2}}|^{\frac{1}{q}-\frac{1}{p}}
\left(\sum_{k\in S_{\frac{n}{2},\frac{n}{2}}}|x_k|^p\right)^{1/p}\right\}\\
\nonumber &= \max\left\{ 1,(n+1)^{\frac1q-\frac1p} 2^\frac{1}{p}\right\}.
\end{align*}
With the choice of $n$ above, we see that
$$
(n+1)^{\frac1q-\frac1p} 2^\frac{1}{p} <1.
$$
Therefore $\|x\|_{\ell^p_q}=1$. Similarly, one may verify that $\|y\|_{\ell^p_q}=1.$
Moreover, we may observe that
$$
\|x+y\|_{\ell^p_q}=2 \quad \text{and} \quad \|x-y\|_{\ell^p_q}=2.
$$
Hence, we obtain
$$
\frac{\|x+y\|_{\ell^p_q}^2+\|x-y\|_{\ell^p_q}^2}{2(\|x\|_{\ell^p_q}^2+\|y\|_{\ell^p_q}^2)}
=\frac{2^2+2^2}{2(1^2+1^2)}=2.
$$
Hence we conclude that
$C_{\rm NJ}(\ell^p_q)= 2$. With the same choices of $x$ and $y$, we have
$$
C_{\rm J}(\ell^p_q)=\sup\{\min \{\|x+y\|_{\ell^p_q},\|x-y\|_{\ell^p_q}\}\colon x,y\in X,
\|x\|_{\ell^p_q}=\|y\|_{\ell^p_q}=1\}=2.
$$

We shall now consider the general case where $d\ge1$.
Let $n\in\mathbb{Z}$ be an even number with
$
n> 2^{\frac{q}{d(q-p)}}-1,
$
or equivalently
$$
(n+1)^{d(\frac1q-\frac1p)}<2^{-\frac1p}.
$$
Let $x\in \ell^p_q$ be the function
$x \colon \mathbb{Z}^d\rightarrow \mathbb{R}$ where
$$
x(k):=\left\{
\begin{array}{ll}
1, &\text{if $k=(0,0,\dots,0), (n,0,\dots,0)$}\\
0, &\text{otherwise}.
\end{array}\right.
$$
and  $y\in \ell^p_q$ be the function
$y \colon \mathbb{Z}^d\rightarrow \mathbb{R}$ where
$$
y(k):=\left\{
\begin{array}{ll}
1, &\text{if $k=(0,0,\dots,0)$}\\
-1, &\text{if $k=(n,0,\dots,0)$}\\
0, &\text{otherwise}.
\end{array}\right.
$$
Then, we have
\begin{align*}
\|x\|_{\ell^p_q}
\nonumber
&=\sup_{m\in\mathbb{Z}^d, N\in \omega}
|S_{m,N}|^{\frac1q-\frac1p} \left(\sum_{k\in S_{m,N}} |x_k|^p\right)^\frac1p\\
&=\max\left\{ 1, |S_{\frac{n}{2},\frac{n}{2}}|^{d(\frac{1}{q}-\frac{1}{p})}
\left(\sum_{k\in S_{\frac{n}{2},\frac{n}{2}}}|x_k|^p\right)^{1/p}\right\}\\
\nonumber &= \max\left\{1, (n+1)^{d(\frac1q-\frac1p)} 2^\frac{1}{p}\right\}.
\end{align*}
Note that with the choice of $n$ above, we have
$$
(n+1)^{d(\frac1q-\frac1p)} 2^\frac{1}{p} <1,
$$
whence $\|x\|_{\ell^p_q}=1$. Similarly $\|y\|_{\ell^p_q}=1.$
Moreover, we also have
$$
\|x+y\|_{\ell^p_q}=2 \quad \text{and} \quad \|x-y\|_{\ell^p_q}=2.
$$
Therefore, we obtain
$$
\frac{\|x+y\|_{\ell^p_q}^2+\|x-y\|_{\ell^p_q}^2}
{2(\|x\|_{\ell^p_q}^2+\|y\|_{\ell^p_q}^2)}=\frac{2^2+2^2}{2(1^2+1^2)}=2,
$$
whence $C_{\rm NJ}(\ell^p_q)=2$. The same choices of $x$ and $y$ give
$$
C_{\rm J}(\ell^p_q)
=\sup\{\min \{\|x+y\|_{\ell^p_q},\|x-y\|_{\ell^p_q}\}\colon x,y\in X,
\|x\|_{\ell^p_q}=\|y\|_{\ell^p_q}=1\}
=2.
$$
Finally as for the Dunkl--Williams constant,
we use the couple
$x+y$ and $(1+r)x+(1-r)y$
for $0<r<1$
and argue similarly to the case of Morrey spaces.
\end{proof}

\bigskip

\noindent{\bf Acknowledgement}. The first author is supported by ITB Research
and Innovation Program 2019. The second author is supported by National Research
Foundation of South Africa, Grant No. 109297. The third author is supported by 
Grant-in-Aid for Scientific Research (C) (Grant No. 16K05209), the Japan Society 
for the Promotion of Science, and by People's Friendship University of Russia.


\begin{thebibliography}{999}
\bibitem{AX12}
D.R.~Adams and J.~Xiao, ``Morrey spaces in harmonic analysis'', \emph{Ark.
Mat.} {\bf 50} (2012), 201--230.

\bibitem{Cassini86}
E.~Casini, ``About some parameters of normed linear spaces'',
\emph{Atti Accad. Naz. Lincei Rend. Cl. Sci. Fis. Mat. Natur. -- Ser. VIII}
{\bf 80} (1986), no.~1-2, 11--15.

\bibitem{DuWi64}
C.F.~Dunkl and K.S.~Williams,
``A simple norm inequality'',
\emph{Amer. Math. Monthly} {\bf 71} (1964), 53--54.

\bibitem{GaLa90}
J.~Gao and K.-S.~Lau, ``On the geometry of spheres in normed linear spaces'',
\emph{J. Aust. Math. Soc. Ser. A} {\bf 48} (1990), 101--112.

\bibitem{GKS18}
H.~Gunawan, E.~Kikianty, and C.~Schwanke, ``Discrete Morrey spaces and their
inclusion properties'', \emph{Math. Nachr.} {\bf 291} (2018), no.~8-9, 1283--1296.

\bibitem{James64}
R.C.~James, ``Uniformly non-square Banach spaces'', \emph{Ann. of Math.}
{\bf 80} (1964), 542--550.

\bibitem{JLM08}
A.~Jim\'{e}nez-Melado, E.~Liorens-Fuster, and E.M.~Mazcu\~{n}\'{a}n-Navarro,
``The Dunkl--Williams constant, convexity, smoothness and normal structure'',
\emph{J. Math. Anal. Appl.} {\bf 342}, Issue 1 (2008), 298--310.

\bibitem{JoNe35}
P.~Jordan and J.~von Neumann, ``On inner products in linear, metric spaces''.
\emph{Ann. Math.} (2) {\bf 36} (1935), 719--723.

\bibitem{KM01}
M.~Kato and L.~Maligandra, ``On James and Jordan-von Neumann constants of
Lorentz sequence spaces'', \emph{J. Math. Anal. Appl.} {\bf 258} (2001),
457--465.

\bibitem{SaSuTa11}
Y.~Sawano, S.~ Sugano, and H.~Tanaka, ``Olsen's inequality and its applications
to Schr\"odinger equations'', in \emph{Harmonic Analysis and Nonlinear Partial
Differential Equations}, 5180, RIMS K$\hat{\rm o}$ky$\hat{\rm u}$roku Bessatsu, B26,
Res. Inst. Math. Sci. (RIMS), Kyoto, 2011.
\end{thebibliography}
\end{document}